\documentclass[12pt]{amsart}

\usepackage[T1]{fontenc}
\usepackage{lmodern,amsthm,amsmath,amssymb,mathrsfs,mathtools,booktabs,longtable,array}
\mathtoolsset{showonlyrefs=true}
\usepackage[a4paper,left=2cm,right=2cm,top=2.5cm,bottom=2.5cm]{geometry}

\usepackage{cite}

\usepackage{microtype}

\usepackage{xcolor}
\usepackage[unicode,pdfencoding=auto,colorlinks=true,linkcolor=blue,citecolor=blue,urlcolor=blue]{hyperref}

\newtheorem{theorem}{Theorem}[section]
\newtheorem{proposition}[theorem]{Proposition}
\newtheorem{lemma}[theorem]{Lemma}
\newtheorem{corollary}[theorem]{Corollary}
\theoremstyle{remark}
\newtheorem{remark}[theorem]{Remark}
\newtheorem{example}[theorem]{Example}
\numberwithin{equation}{section}
\theoremstyle{definition}
\newtheorem{question}[theorem]{Problem}

\newcommand{\F}{\mathbb F}
\newcommand{\Mat}{\mathrm M}
\newcommand{\Id}{I}
\newcommand{\TS}{\mathcal T}
\newcommand{\IC}{\mathcal E}
\DeclareMathOperator{\diag}{diag}

\title[A cubic obstruction to matrix decompositions]{A cubic obstruction to sums of tripotent\\ and square-zero matrices}

\author[D. V. Phuc]{Dang Vo Phuc$^{1,\, *}$}
\address{$^{1}$Department of Mathematics, FPT University, An Phu Thinh New Urban Area, Quy Nhon, Vietnam}
\email{dangphuc150488@gmail.com}
\thanks{$^{*}$Corresponding author}
\thanks{ORCID (D. V. Phuc): \url{https://orcid.org/0000-0002-6885-3996}}

\author[P. V. Danchev]{Peter V. Danchev$^{2}$}
\address{$^{2}$Institute of Mathematics and Informatics, Bulgarian Academy of Sciences, Bulgaria}
\email{danchev@math.bas.bg; pvdanchev@yahoo.com}
\thanks{ORCID (P. V. Danchev): \url{https://orcid.org/0000-0002-2016-2336}}

\date{}

\begin{document}

\begin{abstract}
Let $R$ be a nonzero unital ring of characteristic $2$, and let $\lambda$ be a central unit of $R$. We prove that an element $a\in R$ satisfying $a^3+a^2+\lambda=0$ cannot be expressed as the sum of a tripotent element and a square-zero element. The proof does not require the summands to commute. Over $\F_2$, this gives counterexamples in every order divisible by three, including order twelve, and proves the optimality of the known uniform exponent four in those orders. We also give constructive proofs that every matrix of order four or five over $\F_2$ is a sum of a tripotent and a square-zero matrix, and deduce a complete criterion in order three. Combining the cubic obstruction with the quartic obstruction for idempotent-plus-cube-zero decompositions, we obtain two-way incomparability of the corresponding decomposition classes in every positive order divisible by twelve.
\end{abstract}

\maketitle
\enlargethispage{2pt}

\section{Introduction}

Additive decompositions of matrices into elements satisfying prescribed polynomial identities depend on the ground field and on the degrees of the defining relations. Breaz et al. \cite[Theorem~3]{BCDM} proved that every square matrix over $\F_2$ is a sum of an idempotent and a nilpotent, and that this property characterizes $\F_2$ among fields. A bound on the nilpotency index changes the problem. \v{S}ter \cite{Ster2018} established the bound four, while Shitov \cite[Theorem~1 and its proof]{Shitov} showed that an odd direct sum of companion matrices of $x^4+x^3+1$ cannot be an idempotent plus a cube-zero matrix. \v{S}ter's subsequent ring-theoretic argument \cite[Theorem~3 and Corollary~4]{Ster2022} established the necessity of the bound four in every positive order divisible by four.

Tripotent summands provide a related problem. Abdolyousefi and Chen \cite[Lemma~4.1]{AbdolyousefiChen2018} characterized the fields over which every matrix is a sum of a tripotent and a nilpotent: they are $\F_2$ and $\F_3$. Their result does not require the nilpotent summand to have square zero. For decompositions with a square-zero summand, de Seguins Pazzis \cite[Theorem~1]{Pazzis} obtained a criterion with an idempotent first summand over an arbitrary field. Breaz \cite{Breaz2018} studied periodic and nilpotent summands over finite fields. For a finite field $F$, Danchev et al. \cite[Theorem~3.6]{DGLDiagonalizable} proved that every matrix of order two, and every matrix of order $n\geq3$ with $|F|\geq n+1$, is a sum of a matrix diagonalizable over $F$ and a square-zero matrix. They also proved existence of potent-plus-square-zero decompositions over finite fields and over finite commutative rings with square-zero Jacobson radical, without prescribing a potency exponent independent of the matrix order \cite[Proposition~2.2 and Theorem~2.6]{DGLPotent}.

Further work gives rank criteria for torsion-plus-square-zero decompositions of nilpotent endomorphisms \cite[Theorem~2.2]{DGLRoots2023} and of periodic matrices over fields of positive characteristic \cite[Theorem~2.8]{DGLPeriodic2025}. These results do not prescribe the exponent required for tripotency. Under appropriate block hypotheses, a square-zero perturbation can also prescribe a characteristic polynomial \cite[Theorem~2.3]{DGLPrescribed2024}; a prescribed characteristic polynomial alone need not force the identity $T^3=T$, since the primary-block sizes matter. Square-zero matrices also occur in the analysis of commuting solutions of matrix polynomial equations (see \cite{PhucMordell}).

For finite fields of odd cardinality at least five, Danchev et al. \cite[Theorem~2.2]{DGLOdd} proved that every matrix is a sum of a diagonalizable matrix and a square-zero matrix. Over $\F_3$, they obtained an obstruction from direct sums of companion matrices of a fixed irreducible cubic of nonzero trace \cite[Proposition~3.6]{DGLOdd}. Their recent characteristic-two preprint proves the diagonalizable-plus-square-zero result over finite fields of cardinality at least four and establishes, over $\F_2$, the uniform decomposition $A=D+N$ with $D^4=D$ and $N^2=0$ \cite[Theorem~3.1 and Corollary~3.3]{DGLEven}. The identity $D^4=D$ does not imply $D^3=D$. It therefore remains necessary to determine when a matrix over $\F_2$ admits
\begin{equation}\label{eq:target}
 A=T+N,\qquad T^3=T,\qquad N^2=0.
\end{equation}
In characteristic two, $x^3-x=x(x+1)^2$, so tripotency differs from idempotency. The orders $8n+4$ appearing in Shitov's obstruction suggest asking whether every matrix in $\Mat_{8n+4}(\F_2)$ admits \eqref{eq:target}. In the present paper, we show that this fails at order twelve, whereas order four has an affirmative answer. We now study a cubic obstruction to \eqref{eq:target} in characteristic two and explicit constructions in orders four and five. The main results and their consequences are as follows.
\begin{itemize}
\item[$\bullet$] If $R$ is a nonzero unital ring of characteristic $2$, $\lambda\in Z(R)\cap R^\times$, and $a\in R$ satisfies $a^3+a^2+\lambda=0$, then $a$ is not a sum of a tripotent and a square-zero element (Theorem~\ref{thm:cubic}, Section~\ref{sec:obstruction}).

\item[$\bullet$] Every matrix in $\Mat_d(\F_2)$, with $d\in\{4,5\}$, admits \eqref{eq:target}, by explicit companion-matrix constructions and reductions through primary rational canonical form (Propositions~\ref{prop:order-four} and~\ref{prop:order-five}, Section~\ref{sec:matrices}). The cubic obstruction and the smaller decomposable blocks also yield a complete criterion in order three (Corollary~\ref{cor:order-three}, Section~\ref{sec:obstruction}).

\item[$\bullet$] As consequences of the cubic obstruction, a fixed matrix $C\in\Mat_3(\F_2)$ gives counterexamples $C^{\oplus r}$ in every order $3r$, including the orders $8n+4$ with $n\equiv1\pmod3$. Together with the known four-potent bound, this proves that four is the least uniform potency exponent in order $3r$ (Theorem~\ref{thm:family} and Corollaries~\ref{cor:conjecture} and~\ref{cor:optimal}, Section~\ref{sec:matrices}). Combining the cubic and quartic obstructions with explicit positive decompositions proves that neither of the two decomposition classes contains the other in any positive order divisible by twelve (Theorem~\ref{thm:comparison}, Section~\ref{sec:matrices}).
\end{itemize}

The proof of the cubic obstruction adapts the sandwiched polynomial identities used in characteristic three in \cite[Lemma~3.5 and Proposition~3.6]{DGLOdd}. The characteristic-two cancellation forces $uau=ua^2u=0$ for a hypothetical square-zero summand $u$; an explicit inverse of $a^3+a$ then gives a contradiction. This argument applies to the full ring and therefore allows summands that mix the blocks of a direct sum. For the comparison, we use the quartic obstruction of \v{S}ter \cite[Theorem~3]{Ster2022}. 

The remainder of the paper is organized as follows. Section~\ref{sec:prelim} introduces the decomposition classes and the canonical-form notation, establishes the identities needed for the obstruction argument, and collects the known results and decompositions of small primary blocks used in the applications. Section~\ref{sec:obstruction} develops the argument at the ring level before specializing it to matrix rings. The comparison of the two irreducible cubics over $\F_2$ then leads to the complete criterion in order three. Finally, Section~\ref{sec:matrices} begins with the direct-sum construction and an explicit similarity calculation for the order-six example. These applications are followed by the comparison of the two decomposition classes and the optimality argument for the uniform potency exponent. The constructive part of the section treats companion matrices of orders four and five and extends their decompositions to arbitrary matrices by examining the possible primary configurations. The section concludes with open problems concerning the remaining orders and a general membership criterion in terms of elementary divisors.

\section{Preliminaries}\label{sec:prelim}

In this section, we fix our ring-theoretic and matrix conventions and collect the auxiliary results needed in Sections~\ref{sec:obstruction} and~\ref{sec:matrices}. After introducing the decomposition sets and canonical-form notation, we prove two elementary noncommutative identities in characteristic two: one expands products involving a square-zero element while preserving factor order, and the other provides an explicit inverse for $a^3+a$ when $a^3+a^2+\lambda=0$. We then record the quartic obstruction of \v{S}ter \cite{Ster2022} and the four-potent decomposition over $\F_2$ of Danchev et al. \cite{DGLEven}. The final lemma establishes explicit tripotent-plus-square-zero decompositions for small primary blocks.

All rings are associative and have an identity element $1$. The assumption that a ring is nonzero means $1\neq0$. For a ring $R$, its center is denoted by $Z(R)$ and its group of units by $R^\times$. An element $t\in R$ is tripotent if $t^3=t$, and an element $u\in R$ is square-zero if $u^2=0$. A cube-zero element is an element $v$ satisfying $v^3=0$. No commutativity of summands is included in these definitions. For a positive integer $d$, let $\Id_d$ denote the identity matrix of order $d$, and set
\begin{align}
 \TS_d&=\{T+N:T,N\in\Mat_d(\F_2),\ T^3=T,\ N^2=0\},\label{eq:Tset}\\
 \IC_d&=\{E+Q:E,Q\in\Mat_d(\F_2),\ E^2=E,\ Q^3=0\}.\label{eq:Eset}
\end{align}
The notation $A^{\oplus r}$ denotes the direct sum of $r$ copies of a matrix $A$. Both classes are invariant under similarity and are preserved by direct sums, because conjugation and direct sums preserve the polynomial identities defining the summands.

For a monic polynomial $f\in\F_2[x]$ of positive degree, let $C(f)$ denote its companion matrix, with ones on the subdiagonal. For a matrix $A$, let $m_A(x)$ and $\chi_A(x)$ denote its minimal and characteristic polynomials, respectively. A matrix $A\in\Mat_d(\F_2)$ is cyclic if there is a vector $v$ such that $v,Av,\ldots,A^{d-1}v$ is a basis. In that basis, $A$ is a companion matrix. We shall use the primary rational canonical form, as recalled in \cite[Section~2]{DGLPotent}: every matrix over $\F_2$ is similar to
\begin{equation}\label{eq:primary-form}
 \bigoplus_{i=1}^{s} C(f_i^{e_i}),
\end{equation}
where each $f_i\in\F_2[x]$ is monic and irreducible and each $e_i$ is a positive integer. Each polynomial $f_i^{e_i}$ is an elementary divisor and divides the minimal polynomial. We write $E_{ij}$ for the usual matrix units, so $E_{ij}E_{k\ell}=\delta_{jk}E_{i\ell}$.

\begin{lemma}\label{lem:expansions}
Let $R$ be a ring of characteristic $2$, let $a,u\in R$ satisfy $u^2=0$, and put $t=a+u$. Then
\begin{align}
 u(t^3+t)u&=ua^3u+uauau+uau,\label{eq:sandwich3}\\
 u(t^4+t^2)u&=ua^4u+ua^2uau+uaua^2u+ua^2u,\label{eq:sandwich4}\\
 (t^3+t)u&=(a^3+a)u+auau+ua^2u.\label{eq:right3}
\end{align}
\end{lemma}

\begin{proof}
Distributivity and $u^2=0$ give
\begin{align*}
 t^2&=a^2+au+ua,\\
 t^3&=a^3+a^2u+aua+ua^2+uau.
\end{align*}
Multiplying the second equality on the left and right by $u$, the terms $a^2u$, $ua^2$, and $uau$ vanish because the resulting products contain $u^2$. The remaining terms are $ua^3u$ and $uauau$. Since $utu=uau$, this proves \eqref{eq:sandwich3}.

Multiplying the displayed expression for $t^3$ on the right by $a+u$ gives
\[
\begin{split}
 t^4={}&a^4+a^3u+a^2ua+aua^2+auau\\
       &+ua^3+ua^2u+uaua.
\end{split}
\]
After multiplication on both sides by $u$, only the terms $a^4$, $a^2ua$, and $aua^2$ remain. Also, $ut^2u=ua^2u$. Adding these expressions proves \eqref{eq:sandwich4}.

Finally, right multiplication of the expression for $t^3$ by $u$ gives
\[
 t^3u=a^3u+auau+ua^2u,
\]
because the terms $a^2u^2$ and $uau^2$ vanish. The equality $tu=au$ now yields \eqref{eq:right3}.
\end{proof}

\begin{lemma}\label{lem:inverse}
Let $R$ be a unital ring of characteristic $2$, let $\lambda\in Z(R)\cap R^\times$, and let $a\in R$ satisfy $a^3+a^2+\lambda=0$. Then
\begin{equation}\label{eq:a-powers}
 a^3=a^2+\lambda,\qquad a^4=a^2+\lambda a+\lambda,
\end{equation}
and $a^3+a$ is a unit with inverse
\begin{equation}\label{eq:explicit-inverse}
 (a^3+a)^{-1}=\lambda^{-2}(a^2+\lambda).
\end{equation}
\end{lemma}

\begin{proof}
The first equality in \eqref{eq:a-powers} is the assumed polynomial identity. Multiplying it by $a$ gives $a^4=a^3+\lambda a$, and substituting $a^3=a^2+\lambda$ proves the second equality. Since $\lambda$ is central and the characteristic is $2$, we have
\begin{align*}
 (a^2+a+\lambda)(a^2+\lambda)
 &=a^4+a^3+\lambda a+\lambda^2\\
 &=\lambda^2.
\end{align*}
For the last equality, use $a^4=a^3+\lambda a$. The two factors in the product commute, since each is a polynomial in $a$ with central coefficients. Consequently, multiplication in either order by $\lambda^{-2}(a^2+\lambda)$ gives $1$. The identity $a^3+a=a^2+a+\lambda$ proves \eqref{eq:explicit-inverse}.
\end{proof}

The comparison theorem in Section~\ref{sec:matrices} uses the particular quartic matrix
\begin{equation}\label{eq:H}
 H=\begin{pmatrix}
 0&0&0&1\\
 1&0&0&0\\
 0&1&0&0\\
 0&0&1&1
 \end{pmatrix}\in\Mat_4(\F_2).
\end{equation}
For the quartic obstruction in Theorem~\ref{thm:comparison}, we use the following result of \v{S}ter \cite[Theorem~3, p.~297]{Ster2022}.

\begin{lemma}\label{lem:quartic-obstruction}
Let $R$ be a nonzero unital ring of characteristic $2$, and let $a\in R$ satisfy $a^4+a^3+1=0$. There do not exist $e,q\in R$ such that $a=e+q$, $e^2=e$, and $q^3=0$.
\end{lemma}

The matrix $H$ is $C(x^4+x^3+1)$, so its companion identity is
\begin{equation}\label{eq:H-quartic}
 H^4+H^3+\Id_4=0.
\end{equation}
Consequently Lemma~\ref{lem:quartic-obstruction} applies to every positive direct sum of copies of $H$.

The next result supplies the upper bound in Corollary~\ref{cor:optimal}. It is \cite[Corollary~3.3]{DGLEven}.

\begin{theorem}\label{thm:fourpotent}
Let $d\geq1$ and let $A\in\Mat_d(\F_2)$. There exist $D,N\in\Mat_d(\F_2)$ such that
\[
 A=D+N,\qquad D^4=D,\qquad N^2=0.
\]
\end{theorem}

The following explicit decompositions will be used in Corollary~\ref{cor:order-three} and in the reductions to cyclic matrices in Propositions~\ref{prop:order-four} and~\ref{prop:order-five}.

\begin{lemma}\label{lem:small-blocks}
Let $d\in\{1,2,3\}$ and $A\in\Mat_d(\F_2)$. Suppose that either $d\leq2$, or $d=3$ and $m_A(x)=(x+\varepsilon)^3$ for some $\varepsilon\in\F_2$. Then $A\in\TS_d$.
\end{lemma}

\begin{proof}
A scalar matrix over $\F_2$ is idempotent and hence tripotent, so it has the required decomposition with zero square-zero summand. This proves the assertion for $d=1$. For $d=2$, if every block in \eqref{eq:primary-form} has order one, the resulting diagonal matrix is idempotent, and similarity preserves this identity. Otherwise the primary rational canonical form consists of a single companion matrix of order two. By similarity invariance, it remains to treat that case. For $a,b\in\F_2$, write
\[
 \begin{pmatrix}0&a\\1&b\end{pmatrix}
 =U+V,\qquad
 U=\begin{pmatrix}0&b+1\\1&b\end{pmatrix},\qquad
 V=\begin{pmatrix}0&a+b+1\\0&0\end{pmatrix}.
\]
The matrix $V$ has square zero. If $b=0$, then $U^2=\Id_2$, while if $b=1$, then $U^2=U$. Thus $U^3=U$ in both cases.

For the remaining case, put $J=E_{21}+E_{32}\in\Mat_3(\F_2)$. We have $J^2=E_{31}\neq0$ and $J^3=0$, so $J+\varepsilon\Id_3$ has minimal polynomial $(x+\varepsilon)^3$. The assumed minimal polynomial and \eqref{eq:primary-form} imply that $A$ is similar to $J+\varepsilon\Id_3$. The identities
\[
 J=(E_{12}+E_{21})+(E_{12}+E_{32}),\qquad
 \Id_3+J=(\Id_3+E_{21})+E_{32}
\]
give the required decompositions. Indeed,
\[
 (E_{12}+E_{21})^2=E_{11}+E_{22},\qquad
 (E_{12}+E_{21})^3=E_{12}+E_{21},\qquad
 (E_{12}+E_{32})^2=0,
\]
whereas $(\Id_3+E_{21})^2=\Id_3$ and $E_{32}^2=0$. Hence the first summand in each decomposition is tripotent and the second is square-zero. Similarity invariance completes the proof.
\end{proof}

\section{The cubic obstruction}\label{sec:obstruction}

We now prove the general cubic obstruction in nonzero unital rings of characteristic two using Lemmas~\ref{lem:expansions} and~\ref{lem:inverse}, and specialize it to matrix rings over arbitrary fields of characteristic two. Comparing the two monic irreducible cubics over $\F_2$ then shows why the quadratic term is essential, and combining this distinction with the smaller primary blocks in Lemma~\ref{lem:small-blocks} yields a complete criterion for tripotent-plus-square-zero decompositions in order three.

\begin{theorem}\label{thm:cubic}
Let $R$ be a nonzero unital ring of characteristic $2$, let $\lambda\in Z(R)\cap R^\times$, and let $a\in R$ satisfy
\begin{equation}
 a^3+a^2+\lambda=0.
\end{equation}
There do not exist $t,u\in R$ such that
\[
 a=t+u,\qquad t^3=t,\qquad u^2=0.
\]
\end{theorem}

\begin{proof}
Suppose that such $t$ and $u$ exist. Since $R$ has characteristic $2$, we have $t=a+u$. Put
\begin{equation}\label{eq:XYZ}
 X=uau,\qquad Y=ua^2u,\qquad Z=uauau.
\end{equation}
The equality $t^3=t$ and \eqref{eq:sandwich3} give
\[
 0=ua^3u+Z+X.
\]
By \eqref{eq:a-powers}, centrality of $\lambda$, and $u^2=0$,
\[
 ua^3u=u(a^2+\lambda)u=Y+\lambda u^2=Y.
\]
It follows that
\begin{equation}\label{eq:Yrelation}
 Y=X+Z.
\end{equation}

Multiplying $t^3=t$ by $t$ yields $t^4=t^2$. Formula \eqref{eq:sandwich4} therefore gives
\[
 0=ua^4u+Yau+uaY+Y.
\]
Using the second equality in \eqref{eq:a-powers}, we obtain
\[
 ua^4u=u(a^2+\lambda a+\lambda)u=Y+\lambda X.
\]
The two copies of $Y$ cancel, so
\begin{equation}\label{eq:Xrelation}
 0=\lambda X+Yau+uaY.
\end{equation}
We now compute both mixed products without changing the order of any factors. By \eqref{eq:Yrelation} and \eqref{eq:XYZ},
\begin{align*}
 Yau&=(X+Z)au=uauau+uauauau,\\
 uaY&=ua(X+Z)=uauau+uauauau.
\end{align*}
Thus $Yau=uaY$. Equation \eqref{eq:Xrelation} reduces to $\lambda X=0$, because the characteristic is $2$. Since $\lambda$ is a unit, $X=0$. Consequently,
\[
 Z=Xau=0,
\]
and \eqref{eq:Yrelation} gives $Y=0$. We have proved
\begin{equation}\label{eq:vanishing}
 uau=0,\qquad ua^2u=0.
\end{equation}

Right multiplication of $t^3+t=0$ by $u$ and \eqref{eq:right3} now yield
\begin{align*}
 0&=(a^3+a)u+auau+ua^2u\\
  &=(a^3+a)u+aX+Y\\
  &=(a^3+a)u.
\end{align*}
By Lemma~\ref{lem:inverse}, the element $a^3+a$ is a unit. Multiplying by its inverse gives $u=0$. Hence $a=t$ and $a^3+a=t^3+t=0$. This is impossible because $a^3+a$ is a unit in the nonzero ring $R$.
\end{proof}

Taking $R=\Mat_d(F)$ gives a matrix version with no restriction on the cardinality of the characteristic-two field.

\begin{corollary}\label{cor:field}
Let $F$ be a field of characteristic $2$, let $d\geq1$, and let $\lambda\in F\setminus\{0\}$. If $A\in\Mat_d(F)$ satisfies
\[
 A^3+A^2+\lambda\Id_d=0,
\]
then there do not exist $T,N\in\Mat_d(F)$ such that $A=T+N$, $T^3=T$, and $N^2=0$.
\end{corollary}

\begin{proof}
The ring $\Mat_d(F)$ is nonzero and has characteristic $2$. The scalar matrix $\lambda\Id_d$ is central and invertible. The conclusion follows from Theorem~\ref{thm:cubic}.
\end{proof}

\begin{remark}\label{rem:other-cubic}
The two monic irreducible cubics over $\F_2$ are $x^3+x^2+1$ and $q(x)=x^3+x+1$. Indeed, a monic irreducible cubic has constant coefficient $1$ and does not vanish at $1$. Writing it as $x^3+\alpha x^2+\beta x+1$, these conditions give $\alpha+\beta=1$, which yields exactly the two displayed polynomials; both have no root in $\F_2$ and hence are irreducible.

The companion matrix of $q$ has a different decomposition behavior. Set
\[
 M=\begin{pmatrix}0&0&1\\1&0&1\\0&1&0\end{pmatrix},\qquad
 T_0=\begin{pmatrix}0&0&1\\0&0&1\\0&1&0\end{pmatrix},\qquad
 N_0=\begin{pmatrix}0&0&0\\1&0&0\\0&0&0\end{pmatrix}.
\]
Then $M=C(q)=T_0+N_0$, and direct multiplication gives
\[
 T_0^2=\begin{pmatrix}0&1&0\\0&1&0\\0&0&1\end{pmatrix},
 \qquad T_0^2T_0=T_0,\qquad N_0^2=0.
\]
Consequently $M\in\TS_3$, and taking direct sums gives $M^{\oplus r}\in\TS_{3r}$ for every $r\geq1$. Thus irreducibility and degree three alone do not imply the nonexistence conclusion of Corollary~\ref{cor:field}; the two possible coefficients of $x^2$ in the irreducible cubic distinguish the two examples.
\end{remark}

The two irreducible cubics in Remark~\ref{rem:other-cubic}, together with Lemma~\ref{lem:small-blocks}, account for every primary configuration in order three.

\begin{corollary}\label{cor:order-three}
Let $A\in\Mat_3(\F_2)$. Then $A\in\TS_3$ if and only if $\chi_A(x)\neq x^3+x^2+1$.
\end{corollary}

\begin{proof}
If $\chi_A(x)=x^3+x^2+1$, the Cayley--Hamilton theorem gives $A^3+A^2+\Id_3=0$. Corollary~\ref{cor:field} therefore yields $A\notin\TS_3$.

Conversely, suppose that $\chi_A(x)\neq x^3+x^2+1$, and consider \eqref{eq:primary-form}. If every primary block has order at most two, Lemma~\ref{lem:small-blocks} supplies a decomposition of each block, and direct sums and similarity give $A\in\TS_3$. Otherwise there is one primary block of order three, which occupies the whole matrix. Its elementary divisor is either $(x+\varepsilon)^3$ for some $\varepsilon\in\F_2$, or an irreducible cubic $f$. The first case is covered by Lemma~\ref{lem:small-blocks}. In the second case $f=\chi_A$, so Remark~\ref{rem:other-cubic} and the assumed inequality force $f=x^3+x+1$. The matrix $A$ is then similar to the matrix $M$ of that remark and belongs to $\TS_3$. These alternatives exhaust the primary blocks of total order three.
\end{proof}

\section{Matrix counterexamples and comparison with the quartic obstruction}\label{sec:matrices}

In this section, we specialize Corollary~\ref{cor:field} to $\F_2$ to construct direct-sum counterexamples $C^{\oplus r}$ in every order divisible by three, thereby disproving the universal existence of tripotent-plus-square-zero decompositions in $\Mat_{8n+4}(\F_2)$ whenever $n\equiv1\pmod3$, and we examine a non-block-diagonal example of order six similar to $C\oplus C$. Combining explicit decompositions of $C$ and Shitov's matrix $H$ with the cubic obstruction and Lemma~\ref{lem:quartic-obstruction}, we then show that neither of the decomposition classes $\TS_d$ and $\IC_d$ in \eqref{eq:Tset} and \eqref{eq:Eset} contains the other in any positive order divisible by twelve, and deduce that four is the minimal uniform potency exponent in orders divisible by three. Finally, we prove Propositions~\ref{prop:order-four} and~\ref{prop:order-five} by explicit constructions for cyclic matrices and reductions through primary rational canonical form.

Consider
\begin{equation}\label{eq:C}
 C=\begin{pmatrix}
 0&0&1\\
 0&1&1\\
 1&1&0
 \end{pmatrix}\in\Mat_3(\F_2),
 \qquad p(x)=x^3+x^2+1.
\end{equation}

\begin{theorem}\label{thm:family}
For every positive integer $r$, the matrix $C^{\oplus r}$ does not belong to $\TS_{3r}$.
\end{theorem}

\begin{proof}
Matrix multiplication over $\F_2$ gives
\[
 C^2=\begin{pmatrix}
 1&1&0\\
 1&0&1\\
 0&1&0
 \end{pmatrix},\qquad
 C^3=\begin{pmatrix}
 0&1&0\\
 1&1&1\\
 0&1&1
 \end{pmatrix}.
\]
Adding these matrices gives $C^3+C^2=\Id_3$. For $A=C^{\oplus r}$, multiplication of block diagonal matrices therefore gives
\[
 A^3+A^2=(C^3+C^2)^{\oplus r}=\Id_{3r}.
\]
Corollary~\ref{cor:field}, with $F=\F_2$ and $\lambda=1$, excludes every decomposition of $A$ as a tripotent plus a square-zero matrix.
\end{proof}

\begin{corollary}\label{cor:conjecture}
Let $s\geq0$ be an integer and put $n=3s+1$. Then
\[
 C^{\oplus(8s+4)}\in\Mat_{8n+4}(\F_2)\setminus\TS_{8n+4}.
\]
\end{corollary}

\begin{proof}
The order of the displayed matrix is
\[
 3(8s+4)=24s+12=8(3s+1)+4=8n+4.
\]
Its nonmembership follows from Theorem~\ref{thm:family} with $r=8s+4$.
\end{proof}

In particular, the explicit matrix
\begin{equation}
 A_{12}=\diag(C,C,C,C)
\end{equation}
is a counterexample at $n=1$. Corollary~\ref{cor:conjecture} concerns the subfamily $n\equiv1\pmod3$; it does not assert a classification in the remaining orders. The initial case $n=0$ is affirmative by Proposition~\ref{prop:order-four} below.

The cubic obstruction can also be recognized in a matrix that is not displayed in block diagonal form. Consider
\begin{equation}
 B=\begin{pmatrix}
 0&0&0&0&0&1\\
 0&0&0&0&1&0\\
 0&0&0&1&0&0\\
 0&1&0&0&0&1\\
 1&0&0&1&1&0\\
 0&1&1&0&1&1
 \end{pmatrix}.
\end{equation}
This matrix satisfies $B^3+B^2+\Id_6=0$, so Corollary~\ref{cor:field} excludes a tripotent-plus-square-zero decomposition directly, without a change of basis or an enumeration of possible summands. The following example also identifies its similarity class.

\begin{example}\label{ex:B6}
The matrix $B$ is similar over $\F_2$ to $C\oplus C$. Its minimal polynomial is $p(x)$, its characteristic polynomial is $p(x)^2$, and $B\notin\TS_6$.
\end{example}

\begin{proof}
Direct multiplication over $\F_2$ gives
\[
 B^2=\begin{pmatrix}
 0&1&1&0&1&1\\
 1&0&0&1&1&0\\
 0&1&0&0&0&1\\
 0&1&1&0&0&1\\
 1&1&0&1&1&0\\
 1&1&1&0&1&1
 \end{pmatrix},\qquad
 B^3=B^2+\Id_6.
\]
Thus Corollary~\ref{cor:field} gives $B\notin\TS_6$. To establish the similarity assertion, let
\[
 S=\begin{pmatrix}
 1&1&0&0&1&0\\
 0&1&0&1&1&0\\
 0&0&0&0&1&0\\
 0&0&0&0&1&1\\
 0&1&1&0&1&0\\
 0&1&0&0&1&1
 \end{pmatrix}.
\]
The determinant of this integer matrix is $-1$, so $S$ is invertible over $\F_2$. Direct multiplication gives
\[
 BS=S(C\oplus C)=
 \begin{pmatrix}
 0&1&0&0&1&1\\
 0&1&1&0&1&0\\
 0&0&0&0&1&1\\
 0&0&0&1&0&1\\
 1&0&1&0&1&1\\
 0&1&1&1&0&1
 \end{pmatrix}.
\]
Thus $S^{-1}BS=C\oplus C$.

Expansion of $\det(x\Id_3-C)$ along its first row, in characteristic $2$, gives
\[
 \det(x\Id_3-C)=x\bigl(x(x+1)+1\bigr)+(x+1)=x^3+x^2+1.
\]
The polynomial $p$ is irreducible by Remark~\ref{rem:other-cubic}. The equality $p(C)=0$ established in the proof of Theorem~\ref{thm:family} then implies that the minimal polynomial of $C$ is $p$. A polynomial annihilates $C\oplus C$ precisely when it annihilates $C$, so the minimal polynomial of $B$ is also $p$. The determinant of a block diagonal matrix is the product of the determinants of its diagonal blocks; hence the characteristic polynomial of $B$ is $p^2$.
\end{proof}

To compare the two decomposition conditions, we first record explicit decompositions in the opposite directions for $C$ and $H$.

\begin{proposition}\label{prop:positive}
The matrices $C$ and $H$ satisfy $C\in\IC_3$ and $H\in\TS_4$.
\end{proposition}

\begin{proof}
Define
\[
 E=\begin{pmatrix}
 1&1&1\\
 0&0&0\\
 0&0&0
 \end{pmatrix},\qquad
 Q=\begin{pmatrix}
 1&1&0\\
 0&1&1\\
 1&1&0
 \end{pmatrix}.
\]
Entrywise addition over $\F_2$ gives $C=E+Q$. Multiplication gives
\[
 E^2=E,\qquad
 Q^2=\begin{pmatrix}
 1&0&1\\
 1&0&1\\
 1&0&1
 \end{pmatrix}.
\]
Every row of $Q^2Q$ is the sum of the first and third rows of $Q$, which are equal. Therefore $Q^3=0$, and $C\in\IC_3$.

For the matrix $H$ in \eqref{eq:H}, set
\[
 T=\begin{pmatrix}
 0&0&0&1\\
 0&0&1&0\\
 0&1&0&0\\
 0&0&0&1
 \end{pmatrix},\qquad
 N=\begin{pmatrix}
 0&0&0&0\\
 1&0&1&0\\
 0&0&0&0\\
 0&0&1&0
 \end{pmatrix}.
\]
Then $H=T+N$. The image of $N$ is contained in the span of the second and fourth standard basis vectors, and $N$ vanishes on both of these vectors. Hence $N^2=0$. Furthermore,
\[
 T^2=\begin{pmatrix}
 0&0&0&1\\
 0&1&0&0\\
 0&0&1&0\\
 0&0&0&1
 \end{pmatrix},\qquad
 T^2T=\begin{pmatrix}
 0&0&0&1\\
 0&0&1&0\\
 0&1&0&0\\
 0&0&0&1
 \end{pmatrix}=T.
\]
Thus $T^3=T$, and $H\in\TS_4$.
\end{proof}

The decompositions in Proposition~\ref{prop:positive} supply the positive memberships for the comparison. We are now in a position to prove Theorem~\ref{thm:comparison}.

\begin{theorem}\label{thm:comparison}
Let $m\geq1$ be an integer and put $d=12m$. Then
\[
 C^{\oplus4m}\in\IC_d\setminus\TS_d,
 \qquad
 H^{\oplus3m}\in\TS_d\setminus\IC_d.
\]
In particular, neither $\TS_d$ nor $\IC_d$ is contained in the other.
\end{theorem}

\begin{proof}
Both matrices have order $d$, since $3(4m)=4(3m)=12m$. Taking direct sums of the matrices $E,Q$ from Proposition~\ref{prop:positive} gives
\[
 C^{\oplus4m}=E^{\oplus4m}+Q^{\oplus4m},\qquad
 (E^{\oplus4m})^2=E^{\oplus4m},\qquad
 (Q^{\oplus4m})^3=0.
\]
Thus $C^{\oplus4m}\in\IC_d$, whereas Theorem~\ref{thm:family} gives $C^{\oplus4m}\notin\TS_d$.

Likewise, the matrices $T,N$ from Proposition~\ref{prop:positive} give
\[
 H^{\oplus3m}=T^{\oplus3m}+N^{\oplus3m},\qquad
 (T^{\oplus3m})^3=T^{\oplus3m},\qquad
 (N^{\oplus3m})^2=0.
\]
Hence $H^{\oplus3m}\in\TS_d$. By \eqref{eq:H-quartic},
\[
 (H^{\oplus3m})^4+(H^{\oplus3m})^3+\Id_d
 =(H^4+H^3+\Id_4)^{\oplus3m}=0.
\]
Lemma~\ref{lem:quartic-obstruction}, applied in $\Mat_d(\F_2)$, excludes an idempotent-plus-cube-zero decomposition of this matrix. Therefore $H^{\oplus3m}\notin\IC_d$, proving both failures of containment.
\end{proof}

The obstruction also determines the smallest exponent in the uniform potency statement in orders divisible by three.

\begin{corollary}\label{cor:optimal}
Let $r\geq1$. The least integer $k\geq2$ such that every matrix $A\in\Mat_{3r}(\F_2)$ admits a decomposition $A=D+N$ with $D^k=D$ and $N^2=0$ is $4$.
\end{corollary}

\begin{proof}
Theorem~\ref{thm:fourpotent} proves the assertion with $k=4$. Theorem~\ref{thm:family} excludes $k=3$. If $D^2=D$, then $D^3=D^2D=D$, so a decomposition with $k=2$ would also give one with $k=3$. The same counterexample therefore excludes $k=2$, proving the claim.
\end{proof}

The initial order in the progression $8n+4$ admits a positive result. For that order, the smaller blocks in Lemma~\ref{lem:small-blocks} combine with a uniform construction for companion matrices of order four.

\begin{proposition}\label{prop:order-four}
For every $A\in\Mat_4(\F_2)$, there exist $T,N\in\Mat_4(\F_2)$ such that $A=T+N$, $T^3=T$, and $N^2=0$.
\end{proposition}

\begin{proof}
We first treat cyclic matrices. By similarity invariance, it suffices to consider
\[
 K=\begin{pmatrix}
 0&0&0&a\\
 1&0&0&b\\
 0&1&0&c\\
 0&0&1&d
 \end{pmatrix},\qquad a,b,c,d\in\F_2.
\]
Put
\[
 N_d=\begin{pmatrix}
 d&d&0&0\\
 1&d&0&0\\
 0&0&1&1\\
 0&0&1&1
 \end{pmatrix},\qquad
 T=K+N_d=
 \begin{pmatrix}
 d&d&0&a\\
 0&d&0&b\\
 0&1&1&c+1\\
 0&0&0&d+1
 \end{pmatrix}.
\]
The two diagonal blocks of $N_d$ have squares
\[
 \begin{pmatrix}d&d\\1&d\end{pmatrix}^{\!2}
 =\begin{pmatrix}d^2+d&0\\0&d^2+d\end{pmatrix}=0,
 \qquad
 \begin{pmatrix}1&1\\1&1\end{pmatrix}^{\!2}=0,
\]
because $d^2=d$. Hence $N_d^2=0$. Multiplication gives
\[
 T^2=\begin{pmatrix}
 d&0&0&a+db\\
 0&d&0&b\\
 0&d+1&1&b+d(c+1)\\
 0&0&0&d+1
 \end{pmatrix}.
\]
In the product $T^2T$, the entries in the first three columns agree with those of $T$, using $d^2=d$ and $d(d+1)=0$. The entries in the last column are, respectively,
\begin{align*}
 da+(a+db)(d+1)&=a,\\
 db+b(d+1)&=b,\\
 (d+1)b+(c+1)+(b+d(c+1))(d+1)&=c+1,\\
 (d+1)^2&=d+1.
\end{align*}
Therefore $T^3=T$. Since $K=T+N_d$ in characteristic two, this proves the assertion for every cyclic matrix of order four.

Now apply \eqref{eq:primary-form} to an arbitrary matrix $A$ of order four. If a primary block has order four, then $A$ is cyclic and has already been treated. If a primary block has order three, its elementary divisor is either $(x+\varepsilon)^3$ for some $\varepsilon\in\F_2$ or an irreducible cubic $f$. In the first case, Lemma~\ref{lem:small-blocks} applies to that block and to the remaining scalar block. Taking the direct sum of these two decompositions settles this case.

In the second case, $A$ is similar to $A_0=C(f)\oplus[\varepsilon]$ for some $\varepsilon\in\F_2$. We verify that $A_0$ is cyclic. Let $e_1,e_2,e_3$ be the standard basis of $\F_2^3$, put $v=(e_1,1)$, and set
\[
 W=\operatorname{span}_{\F_2}\{v,A_0v,A_0^2v,A_0^3v\}.
\]
Since $f$ is irreducible, $f(\varepsilon)=1$. The companion identity $f(C(f))=0$ therefore gives $f(A_0)v=(0,1)\in W$. Moreover,
\begin{align*}
 v+f(A_0)v&=(e_1,0),\\
 A_0v+\varepsilon f(A_0)v&=(e_2,0),\\
 A_0^2v+\varepsilon^2 f(A_0)v&=(e_3,0).
\end{align*}
Thus $W$ contains the four standard basis vectors of $\F_2^3\oplus\F_2$, so $W=\F_2^4$ and $A_0$ is cyclic. The preceding construction applies.

If there is no primary block of order three or four, all primary blocks have order at most two and belong to the corresponding class $\TS_d$ by Lemma~\ref{lem:small-blocks}. Taking direct sums of their tripotent and square-zero summands gives the required decomposition of the primary rational canonical form of $A$. Conjugating both summands by the change-of-basis matrix proves the assertion for $A$ itself.
\end{proof}

The order-four result also enters the reduction in the next order. We first construct decompositions for all companion matrices of order five, and then use Proposition~\ref{prop:order-four} to handle the primary configurations containing the exceptional cubic block.

\begin{proposition}\label{prop:order-five}
For every $A\in\Mat_5(\F_2)$, there exist $T,N\in\Mat_5(\F_2)$ such that $A=T+N$, $T^3=T$, and $N^2=0$.
\end{proposition}

\begin{proof}
We first treat cyclic matrices. By similarity invariance, it suffices to consider
\[
 K=\begin{pmatrix}
 0&0&0&0&a\\
 1&0&0&0&b\\
 0&1&0&0&c\\
 0&0&1&0&\delta\\
 0&0&0&1&\varepsilon
 \end{pmatrix},
 \qquad a,b,c,\delta,\varepsilon\in\F_2.
\]
Let $e_1,\ldots,e_5$ denote the standard column vectors of $\F_2^5$.

Suppose first that $\varepsilon=0$, and put
\[
 N=\begin{pmatrix}
 0&1&0&b&a\\
 0&0&0&0&0\\
 0&1&0&1+\delta&c\\
 0&0&0&0&0\\
 0&0&0&0&0
 \end{pmatrix},
 \qquad
 T=K+N=\begin{pmatrix}
 0&1&0&b&0\\
 1&0&0&0&b\\
 0&0&0&1+\delta&0\\
 0&0&1&0&\delta\\
 0&0&0&1&0
 \end{pmatrix}.
\]
The image of $N$ is contained in $\operatorname{span}_{\F_2}\{e_1,e_3\}$, and $Ne_1=Ne_3=0$. Thus $N^2=0$. Set $u=be_1+\delta e_3+e_5$. Since $\delta^2=\delta$ and $\delta(\delta+1)=0$, multiplication gives
\[
\begin{aligned}
 T^2+\Id_5&=
 \begin{pmatrix}
 0&0&b&0&b(\delta+1)\\
 0&0&0&0&0\\
 0&0&\delta&0&0\\
 0&0&0&0&0\\
 0&0&1&0&\delta+1
 \end{pmatrix}
 =u\begin{pmatrix}0&0&1&0&\delta+1\end{pmatrix}.
\end{aligned}
\]
Moreover,
\[
 Tu=bTe_1+\delta Te_3+Te_5
 =be_2+\delta e_4+(be_2+\delta e_4)=0.
\]
Multiplying the displayed factorization on the left by $T$ yields $T^3+T=0$, so $T^3=T$.

Suppose now that $\varepsilon=1$. Set $s=b+c+\delta$ and put
\[
 N=\begin{pmatrix}
 0&0&0&s&s+a\\
 0&1&1&\delta&c\\
 0&1&1&\delta&c\\
 0&0&0&0&0\\
 0&0&0&0&0
 \end{pmatrix},
 \qquad
 T=K+N=\begin{pmatrix}
 0&0&0&s&s\\
 1&1&1&\delta&b+c\\
 0&0&1&\delta&0\\
 0&0&1&0&\delta\\
 0&0&0&1&1
 \end{pmatrix}.
\]
Here the image of $N$ is contained in $\operatorname{span}_{\F_2}\{e_1,e_2+e_3\}$. The first column of $N$ is zero, and its second and third columns coincide. Hence $Ne_1=N(e_2+e_3)=0$, proving $N^2=0$.

Define
\[
 v=e_1+e_2,
 \qquad
 w=se_1+\delta e_2+\delta e_3+e_4+e_5.
\]
Using $\delta^2=\delta$ and $s=b+c+\delta$, multiplication gives
\[
\begin{aligned}
 T^2+\Id_5&=
 \begin{pmatrix}
 1&0&s&s&s(\delta+1)\\
 1&0&\delta&\delta&s+\delta\\
 0&0&\delta&\delta&\delta\\
 0&0&1&1&\delta\\
 0&0&1&1&\delta
 \end{pmatrix}\\
 &=v\begin{pmatrix}1&0&0&0&s\end{pmatrix}
  +w\begin{pmatrix}0&0&1&1&\delta\end{pmatrix}.
\end{aligned}
\]
The first two columns of $T$ coincide, so $Tv=0$. The remaining kernel identity follows from
\begin{align*}
 Tw&=s e_2+\delta e_2+\delta(e_2+e_3+e_4)\\
 &\quad +(s e_1+\delta e_2+\delta e_3+e_5)\\
 &\quad +(s e_1+(b+c)e_2+\delta e_4+e_5)\\
 &=(s+b+c+\delta)e_2=0.
\end{align*}
Multiplying the factorization of $T^2+\Id_5$ on the left by $T$ therefore gives $T^3+T=0$. Thus $T^3=T$ in this case as well. Since $K=T+N$ in characteristic two, every cyclic matrix of order five has the required decomposition.

Now apply \eqref{eq:primary-form} to an arbitrary matrix $A\in\Mat_5(\F_2)$. If a primary block has order five, then $A$ is cyclic and has already been treated. Otherwise every primary block has order at most four. Lemma~\ref{lem:small-blocks}, Corollary~\ref{cor:order-three}, and Proposition~\ref{prop:order-four} show that each such block admits the required decomposition unless its elementary divisor is $p(x)=x^3+x^2+1$. If no block with elementary divisor $p$ occurs, taking direct sums of the decompositions proves the assertion.

Suppose that a block $C(p)$ occurs. It is the only such block, since $A$ has order five. The remaining primary blocks have total order two. If they are two scalar blocks $[\alpha]$ and $[\beta]$, with $\alpha,\beta\in\F_2$, then $C(p)\oplus[\alpha]$ has order four and admits a decomposition by Proposition~\ref{prop:order-four}. The scalar matrix $[\beta]$ is tripotent, so taking the direct sum with its decomposition $[\beta]=[\beta]+[0]$ gives the assertion.

It remains to consider a single primary block $C(h)$ of order two, where $h$ is monic of degree two. Since $p$ is irreducible of degree three, $\gcd(p,h)=1$. We verify that
\[
 A_0=C(p)\oplus C(h)
\]
is cyclic. Let $f_1\in\F_2^3$ and $g_1\in\F_2^2$ be the first standard basis vectors, and put $z=(f_1,g_1)$. Suppose that $q(A_0)z=0$ for a polynomial $q$ of degree at most four. Then
\[
 q(C(p))f_1=0,
 \qquad
 q(C(h))g_1=0.
\]
Division of $q$ by $p$ leaves a remainder of degree at most two. The vectors $f_1,C(p)f_1,C(p)^2f_1$ form the standard basis of $\F_2^3$, so the first equality forces that remainder to be zero. Hence $p\mid q$. The same argument with the basis $g_1,C(h)g_1$ gives $h\mid q$. Coprimality now implies $ph\mid q$. Since $\deg(ph)=5$ and $\deg(q)\leq4$, we obtain $q=0$. Thus the vectors $z,A_0z,\ldots,A_0^4z$ are linearly independent and form a basis of $\F_2^5$. Consequently $A_0$ is cyclic, and the preceding construction applies. Similarity invariance completes the proof.
\end{proof}

\begin{remark}
Lemma~\ref{lem:small-blocks} and Propositions~\ref{prop:order-four} and~\ref{prop:order-five} show that
\[
 \Mat_d(\F_2)=\TS_d
 \qquad\text{for }d\in\{1,2,4,5\},
\]
whereas Theorem~\ref{thm:family} gives
\[
 \Mat_d(\F_2)\neq\TS_d
 \qquad\text{whenever }3\mid d.
\]
Within the progression $8n+4$, Proposition~\ref{prop:order-four} settles the case $n=0$ affirmatively, while Corollary~\ref{cor:conjecture} gives counterexamples whenever $n\equiv1\pmod3$.
\end{remark}

Let $p(x)=x^3+x^2+1$. If $A\in\Mat_d(\F_2)$ satisfies $p(A)=0$, then $m_A$ divides $p$. Since $p$ is irreducible by Remark~\ref{rem:other-cubic} and $m_A$ has positive degree, we have $m_A=p$, and every elementary divisor in \eqref{eq:primary-form} equals $p$. Consequently, $d$ is divisible by $3$. Thus the hypothesis of Corollary~\ref{cor:field} cannot hold over $\F_2$ when $3\nmid d$. These results lead to the following open problems.

\begin{question}\label{q:orders}
Determine all integers $d\geq7$ with $3\nmid d$ for which
\[
 \Mat_d(\F_2)=\TS_d.
\]
In particular, determine which orders $d=8n+4$, with $n\geq2$ and $n\not\equiv1\pmod3$, satisfy this equality.
\end{question}

The first order in the progression specified in Problem~\ref{q:orders} is $20$. The elementary-divisor criterion in Corollary~\ref{cor:order-three} also suggests a corresponding classification in higher orders. Such a classification must account for the interaction between primary summands: although $C\notin\TS_3$ by Theorem~\ref{thm:family}, Proposition~\ref{prop:order-four} gives $C\oplus[0]\in\TS_4$.

\begin{question}\label{q:criterion}
For each integer $d\geq6$, determine necessary and sufficient conditions on the elementary divisors in \eqref{eq:primary-form} for a matrix $A\in\Mat_d(\F_2)$ to belong to $\TS_d$. In particular, determine such conditions when $p$ occurs among the elementary divisors of $A$ and $p(A)\neq0$.
\end{question}


\begin{thebibliography}{99}

\bibitem{AbdolyousefiChen2018}
M.~S.~Abdolyousefi and H.~Chen, \emph{Sums of tripotent and nilpotent matrices}, Bull. Korean Math. Soc. \textbf{55} (2018), 913--920.
DOI: \href{https://doi.org/10.4134/BKMS.b170384}{10.4134/BKMS.b170384}.

\bibitem{Breaz2018}
S.~Breaz, \emph{Matrices over finite fields as sums of periodic and nilpotent elements}, Linear Algebra Appl. \textbf{555} (2018), 92--97.
DOI: \href{https://doi.org/10.1016/j.laa.2018.06.017}{10.1016/j.laa.2018.06.017}.

\bibitem{BCDM}
S.~Breaz, G.~C\u{a}lug\u{a}reanu, P.~Danchev, and T.~Micu, \emph{Nil-clean matrix rings}, Linear Algebra Appl. \textbf{439} (2013), 3115--3119.
DOI: \href{https://doi.org/10.1016/j.laa.2013.08.027}{10.1016/j.laa.2013.08.027}.

\bibitem{DGLDiagonalizable}
P.~Danchev, E.~Garc\'ia, and M.~G.~Lozano, \emph{Decompositions of matrices into diagonalizable and square-zero matrices}, Linear Multilinear Algebra \textbf{70} (2022), 4056--4070.
DOI: \href{https://doi.org/10.1080/03081087.2020.1862742}{10.1080/03081087.2020.1862742}.

\bibitem{DGLPotent}
P.~Danchev, E.~Garc\'ia, and M.~G.~Lozano, \emph{Decompositions of matrices into potent and square-zero matrices}, Internat. J. Algebra Comput. \textbf{32} (2022), 251--263.
DOI: \href{https://doi.org/10.1142/S0218196722500126}{10.1142/S0218196722500126}.

\bibitem{DGLRoots2023}
P.~Danchev, E.~Garc\'ia, and M.~G.~Lozano, \emph{Decompositions of endomorphisms into a sum of roots of the unity and nilpotent endomorphisms of fixed nilpotence}, Linear Algebra Appl. \textbf{676} (2023), 44--55.
DOI: \href{https://doi.org/10.1016/j.laa.2023.07.005}{10.1016/j.laa.2023.07.005}.

\bibitem{DGLPrescribed2024}
P.~Danchev, E.~Garc\'ia, and M.~G.~Lozano, \emph{On prescribed characteristic polynomials}, Linear Algebra Appl. \textbf{702} (2024), 1--18.
DOI: \href{https://doi.org/10.1016/j.laa.2024.08.010}{10.1016/j.laa.2024.08.010}.

\bibitem{DGLPeriodic2025}
P.~Danchev, E.~Garc\'ia, and M.~G.~Lozano, \emph{Decompositions of periodic matrices into a sum of special matrices}, Electron. J. Linear Algebra \textbf{41} (2025), 174--180.
DOI: \href{https://doi.org/10.13001/ela.2025.9099}{10.13001/ela.2025.9099}.

\bibitem{DGLEven}
P.~Danchev, E.~Garc\'ia, and M.~G.~Lozano, \emph{Matrices over finite fields of characteristic $2$ as sums of diagonalizable and square-zero matrices}, preprint (2026), arXiv:2604.15286v1.
DOI: \href{https://doi.org/10.48550/arXiv.2604.15286}{10.48550/arXiv.2604.15286}.

\bibitem{DGLOdd}
P.~Danchev, E.~Garc\'ia, and M.~G.~Lozano, \emph{Matrices over finite fields of odd characteristic as sums of diagonalizable and square-zero matrices}, Linear Algebra Appl. \textbf{730} (2026), 35--50.
DOI: \href{https://doi.org/10.1016/j.laa.2025.10.002}{10.1016/j.laa.2025.10.002}.

\bibitem{Pazzis}
C.~de Seguins Pazzis, \emph{Sums of two triangularizable quadratic matrices over an arbitrary field}, Linear Algebra Appl. \textbf{436} (2012), 3293--3302.
DOI: \href{https://doi.org/10.1016/j.laa.2011.11.026}{10.1016/j.laa.2011.11.026}.

\bibitem{PhucMordell}
\DJ.~V.~Ph\'uc, \emph{The matrix Mordell equation $Y^2+X^3=dI$ over commutative rings: regularity, diagonal lifts, and nilpotent solutions}, Linear Multilinear Algebra \textbf{74} (2026), 1903--1923.
DOI: \href{https://doi.org/10.1080/03081087.2026.2705863}{10.1080/03081087.2026.2705863}.

\bibitem{Shitov}
Y.~Shitov, \emph{The ring $\Mat_{8k+4}(\mathbb Z_2)$ is nil-clean of index four}, Indag. Math. (N.S.) \textbf{30} (2019), 1077--1078.
DOI: \href{https://doi.org/10.1016/j.indag.2019.08.002}{10.1016/j.indag.2019.08.002}.

\bibitem{Ster2018}
J.~\v{S}ter, \emph{On expressing matrices over $\mathbb Z_2$ as the sum of an idempotent and a nilpotent}, Linear Algebra Appl. \textbf{544} (2018), 339--349.
DOI: \href{https://doi.org/10.1016/j.laa.2018.01.015}{10.1016/j.laa.2018.01.015}.

\bibitem{Ster2022}
J.~\v{S}ter, \emph{Nil-clean index of $\Mat_n(\F_2)$}, Linear Algebra Appl. \textbf{632} (2022), 294--307.
DOI: \href{https://doi.org/10.1016/j.laa.2021.10.001}{10.1016/j.laa.2021.10.001}.

\end{thebibliography}
\end{document}